\documentclass[letterpaper, 10 pt, conference]{article}  
\usepackage{amssymb,amsmath,amsfonts,mathrsfs,graphicx,amsthm}

\newtheorem{thm}{Theorem}

\newcommand{\RR}{{\mathbb R}}

\newcommand{\SSS}{{\mathbb  S}}

\newcommand{\abs}[1]{\lvert#1\rvert}

\newcommand{\dv}[2]{{\frac{\partial #1}{\partial #2}}}
\newcommand{\dvv}[2]{{\frac{\partial^2 #1}{\partial {#2}^2}}}
\newcommand{\dotex}{{\tfrac{d}{dt}}}

\newcommand{\VHS}{{V_{\text{\tiny HS}}}}
\newcommand{\GHS}{{\Gamma_{\text{\tiny HS}}}}

\newcommand{\Dh}{{\widehat{D}}}

\title{\LARGE \bf
SO(3)-invariant asymptotic observers for dense depth field estimation based on visual data and known camera motion
}
\author{Nad\`ege Zarrouati
\thanks{Nad\`ege Zarrouati is with DGA, 7-9 rue des Mathurins, 92220 Bagneux, France, PHD candidate at Mines-ParisTech, Centre Automatique et Syst\`emes, Unit\'e Math\'ematiques et Syst\`emes
60, boulevard Saint-Michel
75272 Paris Cedex, France
        {\tt\small nadege.zarrouati@mines-paristech.fr}}%
         \and Emanuel Aldea 
\thanks{Emanuel Aldea is with SYSNAV, 57 rue de Montigny, 27200 Vernon, France
        {\tt\small emanuel.aldea@sysnav.fr }}%
        \and Pierre Rouchon
\thanks{ Pierre Rouchon is with Mines-ParisTech, Centre Automatique et Syst\`emes, Unit\'e Math\'ematiques et Syst\`emes
60, boulevard Saint-Michel
75272 Paris Cedex, France
        {\tt\small pierre.rouchon@mines-paristech.fr}}%
}

\begin{document}

\maketitle
\thispagestyle{empty}
\pagestyle{empty}

\begin{abstract}
In this paper, we use known camera motion associated to a video sequence of a static scene in order to estimate and incrementally refine the surrounding depth field. We exploit the SO(3)-invariance of brightness and depth fields dynamics to customize standard image processing techniques. Inspired by the Horn-Schunck method, we propose a SO(3)-invariant cost to estimate the depth field. At each time step, this provides a diffusion equation on the unit Riemannian sphere that is numerically solved to obtain a real time depth field estimation of the entire field of view. Two asymptotic observers are derived from the governing equations of dynamics, respectively based on optical flow and  depth estimations: implemented on noisy sequences of synthetic images as well as on real data, they perform  a more robust and accurate depth estimation. This approach is  complementary to  most methods employing state observers for range estimation, which uniquely  concern single or isolated feature points.
\end{abstract}

\section{Introduction}

Many vision applications are aimed at assisting in interacting with the environment.
In military as well as in civilian applications, moving in an environment requires topographical knowledge: either in order to avoid obstacles or to engage targets. Since this information is often inaccessible in advance, the real-time computation of a 3D map is a goal that has kept the research community busy for many years.
For example, environment reconstruction is tightly related to the SLAM problem \cite{Smith90}, which is addressed by nonlinear filtering of observed key feature locations (e.g. \cite{EKFMono, Montemerlo03a}), or by bundle adjustment \cite{StrasdatMD10,Strasdat-RSS-10}. However, estimating a sparse point cloud is often insufficient, yet the transition between a discrete local distribution of 3D locations to a continuous depth estimation of the surroundings is an ongoing research topic.
Dynamical systems provide interesting means for incrementally estimating depth information based on the output of vision sensors, since only the current estimates are required, and image batch processing is avoided.

For our work, we are interested in recovering in real-time the depth field around the carrier under the assumptions of known camera motion and known projection model for the onboard monocular camera. The problem of designing an observer to estimate the depth of \textit{single or isolated keypoints} has raised a lot of interest, specifically in the case where the relative motion of the carrier is described by constant known \cite{AbdursulIG04, DahlWLH10}, constant unknown \cite{Heyden09} or time-varying known \cite{chen02, Dixon03, Karagiannis05, deluca08, Astolfi10} affine dynamics.
From a different perspective, the seminal paper of \cite{MatthiesKS89} performs incremental depth field refining for the whole field of view via \emph{iconic} (pixel-wise) Kalman filtering. Video systems, typically found on autonomous vehicles, have successfully used this approach for refining the disparity values obtained by stereo cameras in order to estimate the free space ahead \cite{Hoilund10}. Average optical flow estimations over planar surfaces have also been used for terrain following, in order to stabilize the carrier at a certain pseudo-distance \cite{Herisse10}. Yet, none of these methods provide an accurate dense depth estimation in a general setting concerning the environment and the camera dynamics.

We propose a novel frame of methods relying  on a system of partial differential  equations describing the $SO(3)$-invariant  dynamics of the brightness perceived by the camera and of the depth field of the environment. Based on this invariant kinematic model and the knowledge of the camera motion, these methods provide dense estimations of the depth field at each time-step and exploit such $SO(3)$-invariance.

The present paper is structured as follows. The invariant equations governing the dynamics of the brightness and depth fields are recalled in section~\ref{sec:probstat} and their formulation in pinhole coordinates is given. In section~\ref{sec:VarMeth}, we adapt the Horn-Schunck algorithm  to a variational method providing depth estimation. In section~\ref{sec:Observer}, we propose two asymptotic observers for depth field estimations: the first one is based on standard optical flow measures and the second one enables the refinement of rough or inaccurate depth estimations; we prove their convergence  under geometric  assumptions concerning the camera dynamics and the environment. In section~\ref{sec:implement}, we test these methods on  synthetic data and compare their accuracy, their robustness to noise and their convergence rate; tested on real data, this approach gives promising results.

\section{The $SO(3)$-invariant model }
\label{sec:probstat}

\subsection{The partial differential system on $\SSS^2$ }
\label{subsec:PDS}
The model is based on geometric  assumptions introduced  in  \cite{Bonnabel-Rouchon2009}. We consider a spherical camera, whose motion is known. Linear and angular velocities $v(t)$ and $\omega (t)$ are  expressed in the camera frame. Position of the optical center  in the reference frame $\cal R$ is denoted by $C(t)$. Orientation versus $\cal R$  is given by the quaternion $q(t)$: any vector $\varsigma$ in the camera frame corresponds to the vector $q\varsigma q^*$ in the reference frame $\cal R$ using the identification of  vectors as    imaginary quaternions. We have  thus: $\dotex q = \frac{1}{2} q\omega$. A pixel is labeled by  the unit vector $\eta$  in the camera frame: $\eta$ belongs to the sphere  $\mathbb  S^2$ and  receives the brightness $y(t,\eta)$. Thus at each time $t$, the image produced by the camera is described by the  scalar field $\mathbb S^2\ni \eta \mapsto y(t,\eta)\in\RR$.

The scene is modeled as a closed, $C^1$ and convex  surface $\Sigma$  of $\RR^3$,  diffeomorphic to  $\mathbb  S^2$.   The camera is inside the domain $\Omega\subset \RR^3$ delimited by $\Sigma=\partial\Omega$.   To a point $M\in\Sigma$ corresponds  one and only one camera pixel: if the points of $\Sigma$ are labeled by $s\in\SSS^2$,   for each time $t$, a continuous and invertible transformation $\mathbb  S^2\ni s\mapsto \phi(t,s)\in \mathbb  S^2$ enables to express  $\eta$ as a function of $s$: $\eta=\phi(t,s)$.

\begin{figure}[thpb]
      \centering
      \includegraphics[width=0.6\textwidth]{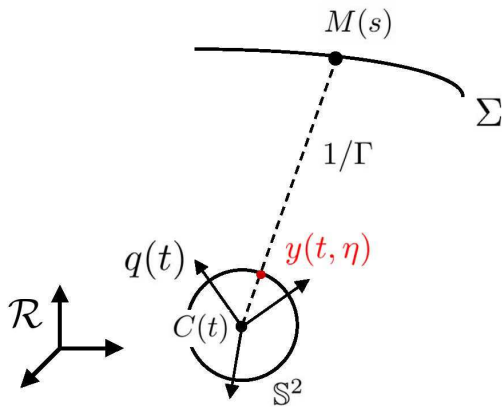}
      \caption{Model and notations of a spherical camera in a static environment.}
      \label{fig:notations}
\end{figure}

The density of light emitted by a point  $M(s)\in\Sigma$ does not depend on the direction of emission ($\Sigma$ is a Lambertian surface) and is independent of $t$ (the scene is static). This means that $y(t,\eta)$ depends only on $s$: thus $y$ can be seen either as a function of $(t,\eta)$ or, via the transformation $\phi$,  as a function of $s$.
 The  distance $C(t)M(s)$   between the optical center and the object seen in the direction $\eta=\phi(t,s)$ is  denoted by $D(t,\eta)$, and its inverse by $\Gamma = 1/D$. Fig.\ref{fig:notations} illustrates the model and the notations. We assume  that $s\mapsto y(s)$ is a $C^1$ function. For each $t$, $s\mapsto D(t,s)$ is $C^1$ since $\Sigma$ is a $C^1$ surface of $\RR^3$.

Under these assumptions, we first have:
\begin{equation}
\left.\dv{y}{t}\right|_s=0, \quad
\left.\dv{\Gamma}{t}\right|_s=\Gamma^2v\cdot\eta
\label{eq:gamma_direction}
\end{equation}
then
\begin{multline*}
\left.\dv{h}{t}\right|_s=\left.\dv{h}{t}\right|_{\eta}+\left.\dv{h}{\eta}\right|_{t}\left.\dv{\eta}{t}\right|_s
=\left.\dv{h}{t}\right|_{\eta}+\nabla h\cdot\left.\dv{\eta}{t}\right|_s
\end{multline*}
where $h$ is any scalar field defined on $\mathbb  S^2$ and $\nabla h$ its gradient with respect to the Riemannian  metric on $\mathbb S^2$. The value of $\nabla h$ at $\eta\in\SSS^2$ is identified with a vector of $\RR^3$ tangent to the sphere at the point $\eta$ also identified  to  a  unitary vector of~$\RR^3$  in the  camera moving frame. The Euclidean scalar product of two vectors $a$ and $b$ in $\RR^3$ is denoted by $a\cdot b$ and their wedge product by $a\times b$. By differentiation, the identity $q\eta q^*=\frac{\overrightarrow{C(t)M(s)}}{\left||C(t)M(s)\right||}$, where $^*$ denotes  conjugation and $\eta$ is identified to an imaginary quaternion,  yields
\begin{equation}
\left.\dv{\eta}{t}\right|_s=\eta \times (\omega +\Gamma\eta\times v)
\end{equation}
since the vector  $\eta\times \omega$ corresponds to the imaginary quaternion $(\omega\eta - \eta\omega)/2$.
 Therefore, the intensity $y(t,\eta)$ and the inverse depth $\Gamma(t,\eta)$ satisfy the following equations:
\begin{equation}
\dv{y}{t} =- \nabla y \cdot(\eta \times (\omega +\Gamma\eta\times v))
\label{eq:flot}
\end{equation}
\begin{equation}
\dv{ \Gamma}{t} =- \nabla \Gamma \cdot(\eta \times (\omega +\Gamma\eta\times v)) + \Gamma^{2}v\cdot\eta
\label{eq:prof}
\end{equation}
Equations \eqref{eq:flot} and \eqref{eq:prof} are $SO(3)$-invariant: they remain unchanged by  any rotation described by the quaternion $\sigma$ and changing $(\eta,\omega,v)$ to
$(\sigma \eta \sigma^*,\sigma \omega \sigma^*,\sigma v \sigma^* )$.
Equation \eqref{eq:flot} is the well-known optical flow equation  that can be found under different forms in numerous papers (see \cite{Murray09} or \cite{Astolfi10} for example), while \eqref{eq:prof} is  less standard (see e.g., \cite{Bonnabel-Rouchon2009}).

\subsection{The system  in pinhole coordinates}
\label{subsec:Pinhole}
To use this model  with camera data, one needs to write  the invariant  equations~\eqref{eq:flot} and~\eqref{eq:prof} with local coordinates on $\SSS^2$ corresponding to a rectangular grid of pixels. One popular solution is to use the pinhole camera model, where the pixel of coordinates $(z_1,z_2)$ corresponds to the unit vector $\eta\in\SSS^2$ of coordinates in $\RR^3$:  $\left(1+z_{1}^2+z_{2}^2\right)^{-1/2}(z_1,z_2,1)^T$. The optical camera axis (pixel $(z_1,z_2)=(0,0)$) corresponds here  to the direction $z_3$. Directions $1$ and $2$ correspond respectively to the horizontal axis from left to right and to the vertical axis from top to bottom on the image frame.

The gradients $\nabla y$ and $\nabla \Gamma$ must be expressed with respect to $z_1$ and $z_2$. Let us detail this derivation for $y$. Firstly, $\nabla y$ is tangent to $\mathbb S^2$, thus $\nabla y\cdot\eta=0$. Secondly, the differential $dy$ corresponds to $\nabla y\cdot d\eta$ and to $\dv{y}{z_{1}} dz_{1} +\dv{y}{z_{2}} dz_{2}$. By identification, we get the Cartesian coordinates of $\nabla y$ in $\RR^3$. Similarly we get the three coordinates of $\nabla\Gamma$. Injecting these expressions in~\eqref{eq:flot} and \eqref{eq:prof}, we get the following partial differential equations (PDE) corresponding to~\eqref{eq:flot} and~\eqref{eq:prof} in local pinhole coordinates:
\begin{equation}
\begin{aligned}
\dv{y}{t}= &- \dv{y}{z_1} \left[\begin{array}{c}z_1z_2\omega_1-(1+{z_1}^2)\omega_2+z_2\omega_3\nonumber\\
+\Gamma\sqrt{1+z_{1}^2+z_{2}^2}(-v_1+z_1v_3)\end{array}\right] \nonumber\\
           &- \dv{y}{z_2} \left[\begin{array}{c}(1+{z_2}^2)\omega_1-z_1z_2\omega_2-z_1\omega_3\\
+\Gamma\sqrt{1+z_{1}^2+z_{2}^2}(-v_2+z_2v_3)\end{array}\right]\nonumber
 \label{eq:flotcart}
\end{aligned}
\end{equation}

\begin{equation}
\begin{aligned}
\dv{\Gamma}{t}= &- \dv{\Gamma}{z_1} \left[\begin{array}{c}z_1z_2\omega_1-(1+{z_1}^2)\omega_2+z_2\omega_3\nonumber\\
+\Gamma\sqrt{1+z_{1}^2+z_{2}^2}(-v_1+z_1v_3)\end{array}\right]\nonumber\\
           &- \dv{\Gamma}{z_2} \left[\begin{array}{c}(1+{z_2}^2)\omega_1-z_1z_2\omega_2-z_1\omega_3\\
+\Gamma\sqrt{1+z_{1}^2+z_{2}^2}(-v_2+z_2v_3)\end{array}\right]\nonumber\\
&+\Gamma^2(z_1v_1+z_2v_2+v_3)\nonumber
 \label{eq:profcart}
\end{aligned}
\end{equation}
where $v_1$, $v_2$, $v_3$, $\omega_1$, $\omega_2$, $\omega_3$ are the components of linear and angular velocities in the camera frame.

\section{Depth estimation inspired by Horn-Schunck method}
\label{sec:VarMeth}

\subsection{The Horn-Schunck variational method}
\label{ssec:HS}
In~\cite{HS81}, Horn and Schunck described a method to compute the optical flow, defined as "the distribution of apparent velocities of movement of brightness patterns in an image". The entire method is based on the optical flow constraint written in  the compact form
 \begin{equation}
\dv{y}{t}+V_1\dv{y}{z_1}+V_2\dv{y}{z_2}=0
\label{eq:OFC}
\end{equation}
Identification with~\eqref{eq:flotcart} yields
\begin{equation}
\begin{split}
&V_1(t,z)=f_{1}(z,\omega(t))+\Gamma(t,z) g_{1}(z,v(t))\\
&V_2(t,z)=f_{2}(z,\omega(t))+\Gamma(t,z) g_{2}(z,v(t))
\end{split}
\label{eq:defV}
\end{equation}
with
\begin{equation}
\begin{split}
&f_{1}(z,\omega)=z_1z_2\omega_1-(1+{z_1}^2)\omega_2+z_2\omega_3\\
&g_{1}(z,v)=\sqrt{1+z_{1}^2+z_{2}^2}(-v_1+z_1v_3)\\
&f_{2}(z,\omega)=(1+{z_2}^2)\omega_1-z_1z_2\omega_2-z_1\omega_3\\
&g_{2}(z,v)=\sqrt{1+z_{1}^2+z_{2}^2}(-v_2+z_2v_3)
.
\end{split}
\label{eq:fg}
\end{equation}
For each time $t$, the apparent velocity  field $V=V_1\dv{}{z_1}+V_2\dv{}{z_2}$ is then estimated  by   minimizing versus $W=W_1\dv{}{z_1}+W_2\dv{}{z_2}$ the following cost (the image $\mathcal I$ is a rectangle of $\RR^2$ here)
\begin{multline}
I(W)=\iint_{\mathcal I}\bigg(\big(\dv{y}{t}+W_1\dv{y}{z_1}+W_2\dv{y}{z_2}\big)^2 \\
+\alpha^{2} ({\nabla W_1}^2+{\nabla W_2}^2)\bigg) dz_1dz_2
\label{eq:integraleOFC1}
\end{multline}
where $\nabla$ is the gradient operator in the Euclidian plane $(z_1,z_2)$,  $\alpha >0$ is a regularization parameter and the
 partial derivatives $\dv{y}{t}$, $\dv{y}{z_1}$ and  $\dv{y}{z_2}$  are assumed to be known.

 Such  Horn-Schunk estimation of  $V$ at time $t$  is denoted by
 $$\VHS(t,z)=\VHS_1(t,z) \dv{}{z_1}+\VHS_2(t,z) \dv{}{z_2}.$$
  For each time $t$, usual  calculus of variation yields the following PDE's for $\VHS$:
\begin{align*}
&\left(\dv{y}{z_1}\right)^{2} \VHS_1+\dv{y}{z_1}\dv{y}{z_2} \VHS_2=\alpha^{2} \Delta \VHS_1-\dv{y}{z_1}\dv{y}{t}
\\
&\dv{y}{z_1}\dv{y}{z_2}\VHS_1 + \left(\dv{y}{z_2}\right) ^{2} \VHS_2=\alpha^{2} \Delta \VHS_2-\dv{y}{z_2}\dv{y}{t}
\end{align*}
with  boundary conditions $\dv{\VHS_1}{n}=\dv{\VHS_2}{n}=0$  ($n$ the normal to $\partial \mathcal I$).
Here  $\Delta$ is the Laplacian operator in the Euclidian space $(z_1,z_2)$.
The numerical resolution is usually based on
\begin{itemize}
 \item computations of $\dv{y}{z_1}$, $\dv{y}{z_2}$ and $\dv{y}{t}$ via  differentiation  filters (Sobel filtering) directly from the image data at different times around $t$.
\item approximation of $\Delta \VHS_1$ and $\Delta \VHS_2$ by the difference between the weighted mean $\bar \VHS_1$ and $\bar \VHS_2$ of $\VHS_1$ and $\VHS_2$  on the neighboring pixels and their values at the  current pixel;
\item iterative resolution (Jacobi scheme) of the resulting linear system in $\VHS_1$ and $\VHS_2$.
\end{itemize}
The convergence of this numerical method of resolution was proven in \cite{convHS}.
Three parameters have a direct impact on the speed of convergence and on the precision: the regularization parameter $\alpha$, the number of iterations for the  Jacobi scheme   and the initial values of $\VHS_1$ and $\VHS_2$ at the beginning of this iteration step. To be specific, $\alpha$ should neither be too small in order to filter noise appearing  in differentiation  filters applied on  $y$, nor too large  in order to have $\VHS$ close to $V$ when $\nabla V\neq 0$.

\subsection{Adaptation to depth estimation }
\label{ssec:HSDepth}

Instead of minimizing the cost $I$ given by~\eqref{eq:integraleOFC1} with respect to any $W_1$ and $W_2$, let us define a new invariant cost  $J$,
\begin{multline}
J(\Upsilon)=
\iint_{\mathcal J} \bigg(\left(\dv{y}{t} +\nabla y \cdot\left(\eta \times \left(\omega +\Upsilon\eta\times v\right)\right)\right)^2\\
+\alpha^{2} {\nabla \Upsilon}^2 \bigg) d\sigma_{\eta}
\label{eq:integraleOFC2}
\end{multline}
and minimize it with respect to any depth profile $\mathcal J \ni \eta \mapsto \Upsilon(t,\eta)\in\RR$. The time $t$ is fixed here and $d\sigma_\eta$ is the Riemannian infinitesimal surface element on $\SSS^2$. $\mathcal J\subset \SSS^2$ is the domain where $y$ is measured and $\alpha >0$ is the  regularization parameter.

The first order stationary condition of $J$ with respect to any variation of $\Upsilon$ yields the following invariant PDE characterizing the resulting  estimation $\GHS$ of   $\Gamma$:
\begin{multline}
\alpha^{2}\Delta \GHS=
\left(\dv{y}{t} +\nabla y \cdot\left(\eta \times \left(\omega +\GHS\eta\times v\right)\right)\right) \ldots \\ \ldots \left(\nabla y \cdot\left(\eta \times \left(\eta\times v\right)\right)\right) \phantom{eeeeee} \text{on }\mathcal J
\label{eq:CV}
\end{multline}
\begin{equation}
\phantom{eeeeeeeeeeeeeeeeeee}\dv{\GHS}{n}=0 \phantom{eeeeeeeeeeeee} \text{on }\partial\mathcal J
\label{eq:bord}
\end{equation}
where $\Delta \GHS$ is the Laplacian of $\GHS$  on the Riemannian sphere $\mathbb  S^2$ and $\partial\mathcal J$ is the  boundary of $\mathcal J$,  assumed to be piece-wise smooth and with unit normal vector $n$.

In pinhole coordinates $(z_1,z_2)$, we have
\begin{align*}
&d\sigma_{\eta}=\left(1+z_{1}^2+z_{2}^2\right)^{-3/2}dz_{1}dz_{2}
\\
&\nabla\Upsilon^2= \left(1+z_{1}^2+z_{2}^2\right)\left(\dv{\Upsilon}{z_{1}}^2+\dv{\Upsilon}{z_{2}}^2+(z_1\dv{\Upsilon}{z_{1}}+z_2\dv{\Upsilon}{z_{2}})^2\right)\\
&\left(\dv{y}{t} +\nabla y \cdot\left(\eta \times \left(\omega +\Upsilon\eta\times v\right)\right)\right)^2
=(F+\Upsilon G)^2
\end{align*}

where
\begin{equation}
\begin{split}
&F=\dv{y}{t}+f_{1}(z,\omega)\dv{y}{z_1}+f_{2}(z,\omega)\dv{y}{z_2}\\
&G=g_{1}(z,v)\dv{y}{z_1}+g_{2}(z,v)\dv{y}{z_2}.
\label{eq:F1F2}
\end{split}
\end{equation}
Consequently, the
first order stationary condition~\eqref{eq:CV}  reads in $(z_1,z_2)$ coordinates:
\begin{equation}
\begin{split}
\GHS G^2+ FG=&\alpha^2\left[\frac{\partial}{z_1}\left(\frac{1+z_1^2}{\sqrt{1+z_{1}^2+z_{2}^2}}\dv{\GHS}{z_1}\right)
\right.\\
&+\frac{\partial}{z_2}\left(\frac{1+z_2^2}{\sqrt{1+z_{1}^2+z_{2}^2}}\dv{\GHS}{z_2}\right)\\
&+\frac{\partial}{z_2}\left(\frac{z_1z_2}{\sqrt{1+z_{1}^2+z_{2}^2}}\dv{\GHS}{z_1}\right)\\
&\left. +\frac{\partial}{z_1}\left(\frac{z_1z_2}{\sqrt{1+z_{1}^2+z_{2}^2}}\dv{\GHS}{z_2}\right)\right]
\end{split}
\label{eq:intrinsicEDP}
\end{equation}
on the rectangular domain $\mathcal I=[-\bar z_1,\bar z_1]\times [-\bar z_2,\bar z_2]$ ($\bar z_1,\bar z_2>0$ with $\bar z_1^2+\bar z_2^2 < 1$).

The right term of~\eqref{eq:intrinsicEDP} corresponds to the Laplacian operator on the Riemannian sphere $\SSS^2$ in pinhole coordinates. The numerical resolution of this scalar diffusion providing the estimation $\GHS$ of $\Gamma$ is similar to the one used for the Horn-Schunck estimation $\VHS$ of $V$. The functional $I(W)$ defined in \eqref{eq:integraleOFC1} is minimized with respect to two varying parameters $W_1$ and $W_2$ while there is really only one unknown function in this problem: the depth field $\Gamma$. On the contrary, the functional $J(\Upsilon)$ takes full advantage of the knowledge of the camera dynamics since the only varying parameter here is $\Upsilon$.

\section{Depth estimation via asymptotic observers}
\label{sec:Observer}

\subsection{Asymptotic observer based on optical flow measures ($\VHS$)}
\label{ssec:ObsVHS}
From any optical flow estimation, such as $\VHS$, it is reasonable to  assume that we have access for each time $t$, to the components in pinhole coordinates of the vector  field
\begin{align}\label{eq:varpi}
\varpi_t:~ \SSS^2\ni \eta \mapsto \varpi_t(\eta)= \eta \times (\omega +\Gamma\eta\times v)\in T_\eta \SSS^2
\end{align}
 appearing in~\eqref{eq:prof}. This vector field can be  considered as a measured output for~\eqref{eq:prof}, expressed as $\varpi_t(\eta)=f_t(\eta)+\Gamma(t,\eta) g_t(\eta)$, where $f_t$ and $g_t$ are the vector fields
\begin{align}
 &f_t:~ \SSS^2\ni \eta \mapsto f_t(\eta)= \eta \times \omega \in T_\eta \SSS^2\label{eq:f}
 \\
 &g_t:~ \SSS^2\ni \eta \mapsto g_t(\eta)= \eta \times (\eta\times v) \in T_\eta \SSS^2\label{eq:g}.
\end{align}
This enables us to propose  the following asymptotic  observer for $D=1/\Gamma$ since it obeys to  $\dv{ D}{t} =- \nabla D \cdot\varpi_t - v\cdot\eta$:
\begin{equation}\label{eq:obserIntrinsic}
\dv{ \Dh}{t} =- \nabla \Dh \cdot\varpi_t - v\cdot\eta + k  g_t\cdot( \Dh f_t+  g_t -\Dh \varpi_t)
\end{equation}
where  $\varpi_t$, $f_t$ and $g_t$ are known time-varying vector fields on $\SSS^2$ and $k>0$ is a tuning parameter. This observer is trivially  $SO(3)$ invariant and reads in  pinhole coordinates:
\begin{equation}
\begin{split}
&\dv{\Dh}{t}=-\dv{\Dh}{z_1}V_1-\dv{\Dh}{z_2}V_2-(z_1v_1+z_2v_2+v_3)\\
&\quad+k  \big(g_1(\Dh f_1+g_1-\Dh V_1)+g_2(\Dh f_2+ g_2-\Dh V_2) \big)
\end{split}
\label{eq:observerPH}
\end{equation}
where $V$ is given by any optical flow estimation and $(f_1,f_2,g_1,g_2)$ are defined  by~\eqref{eq:fg}.

As assumed in the first paragraphs of subsection~\ref{subsec:PDS}, for each time $t$,   there is a one to one smooth mapping between $\eta\in\SSS^2$ attached to the camera pixel and the scene point $M(s)$ corresponding to this pixel.  This means that, for any $t\geq 0$,  the flow $\phi(t,s)$ defined by
\begin{equation}\label{eq:phi}
    \left.\dv{\phi}{t}\right|_{(t,s)} = \varpi_t(\phi(t,s)),\quad \phi(0,s)=s\in\SSS^2
\end{equation}
defines a time varying diffeomorphism on $\SSS^2$. Let us denote by $\phi^{-1}$ the inverse diffeomorphism:
$\phi(t,\phi^{-1}(t,\eta))\equiv \eta$. Assume that $\Gamma(t,\eta) >0 $, $v(t)$ and $\omega(t)$ are   uniformly bounded  for $t\geq 0$ and $\eta\in\SSS^2$. This means that  the trajectory of the camera center $C(t)$ remains strictly inside the convex surface $\Sigma$  with minimal distance to $\Sigma$.  These considerations motivate   the assumptions used in the following theorem.

\begin{thm}\label{thm:error}

Consider $\Gamma(t,\eta)$ associated to the motion of the camera inside the domain $\Omega$ delimited by the scene $\Sigma$, a $C^1$, convex and closed surface as explained in sub-section~\ref{subsec:PDS}. Assume that exist $\bar v>0$, $\bar\omega>0$, $\bar\gamma>0$ and $\bar\Gamma>0$  such that
$$
\forall t\geq 0,~\forall \eta\in\SSS^2,~|v(t)|\leq\bar v,~|\omega(t)|\leq\bar \omega,~\bar\gamma \leq \Gamma(t,\eta) \leq \bar\Gamma.
$$
Then, for $t\geq 0$, $\Gamma(t,\eta)$ is a $C^1$ solution of~\eqref{eq:prof}.
Consider the  observer~\eqref{eq:obserIntrinsic} with  a $C^1$  initial condition  versus $\eta$, $\Dh(0,\eta)$. Then we have the following implications:
\begin{itemize}

\item $\forall t\geq 0$, the solution $\Dh(t,\eta)$ of~\eqref{eq:obserIntrinsic}  exists, is unique and remains $C^1$ versus $\eta$. Moreover
    $$
    t\mapsto \| \Dh(t,\rule{0.2em}{0.2em})- D(t,\rule{0.2em}{0.2em}) \|_{L^\infty} = \max_{\eta\in\SSS^2} \left|\Dh(t,\eta) -D(t,\eta) \right|
    $$
is  decreasing ($L^\infty$ stability).

\item if  additionally for all $s\in\SSS^2$, $\int_0^{+\infty } \|g_\tau(\phi(\tau,s))\|^2 d\tau =+\infty$, then we have
    for all $p >0$,
    $$
    \lim_{t\mapsto +\infty} \int_{\SSS^2} \big|\Dh(t,\eta)-D(t,\eta)\big|^p d\sigma_\eta = 0
    $$
(convergence in any $L^p$ topology)

\item if additionally  there is $\lambda>0$ and $T>0$  such that, for all $t\geq T$ and $s\in\SSS^2$,  $\int_0^{t} \|g_\tau(\phi(\tau,s))\|^2 d\tau \geq \lambda t  $, then we have,
    for all $t\geq T$,
    $$
    \| \Dh(t,\rule{0.2em}{0.2em})- D(t,\rule{0.2em}{0.2em}) \|_{L^\infty} \leq
    e^{-k\bar\gamma\lambda t}
    \| \Dh(0,\rule{0.2em}{0.2em})- D(0,\rule{0.2em}{0.2em}) \|_{L^\infty}
    $$
(exponential convergence in $L^\infty$ topology).

\end{itemize}

\end{thm}
Assumptions on $\int \|g_t(\phi(t,s))\|^2 dt $ can be seen as a condition of persistent excitations. It should be  satisfied for generic motions of the camera.

\begin{proof}
The facts that $v$, $\omega$ and $\Gamma$ are bounded  and that the scene surface  $\Sigma$ is $C^1$, closed and convex,  ensure that the mapping  $\eta=\phi(t,s)$  and its inverse $s=\phi^{-1}(t,\eta)$  are  $C^1$ diffeomorphism on $\SSS^2$ with bounded derivatives versus $s$ and $\eta$ for all time $t>0$. Therefore,  $\Gamma$ is also a function of $(t,s)$. Set $\overline{\Gamma}(t,s)=\Gamma(t,\phi(t,s))$: in the $(t,s)$ independent variables the partial differential equation~\eqref{eq:prof} becomes a set of ordinary differential equations indexed by $s$:
$\left.\dv{\overline{\Gamma}}{t}\right|_s=\overline{\Gamma}^2v(t)\cdot\phi(t,s)$
that reads also $ \left.\dv{(\overline{D})}{t}\right|_s=-v(t)\cdot\phi(t,s)$ with  $\overline{D}=1/\overline{\Gamma}$.
Thus $\overline{D}(t,s) - \overline{D}(0,s) = -\int_{0}^{t} v(\tau )\cdot\phi(\tau,s)~d\tau$. Consequently, $\overline{\Gamma}$ is $C^1$ versus $s$ and thus $\Gamma$ is $C^1$ versus $\eta$.

Set $\overline{\Dh}(t,s)=\Dh(t,\phi(t,s))$. Then
\begin{multline*}
\left.\dv{\overline{\Dh}}{t}\right|_s=- v(t)\cdot\phi(t,s)\\+k \|g_t(\phi(t,s))\|^2\overline{\Gamma}(t,s)(\overline{D}(t,s)-\overline{\Dh}(t,s)).
\end{multline*}
Set $E=\Dh -D $ and $\overline{E}= \overline{\Dh}-\overline{D}$. Then
\begin{equation}\label{eq:E}
\left.\dv{\overline{E}}{t}\right|_s = - k \overline{\Gamma}(t,s) \|g_t(\phi(t,s))\|^2 \overline{E}(t,s)
\end{equation}
Consequently, $\bar E$ is well defined for any $t>0$ and $C^1$ versus $s$. Thus $E$ and consequently $\Dh=E+D$ are also well defined for all $t>0$ and are $C^1$ versus $\eta$.

Since for any $s$ and $t_2> t_1 \geq 0$ we have
$|\overline{E}(t_2,s)| \leq |\overline{E}(t_1,s)|$, we have also
$$
|\overline{E}(t_2,s)|\leq \max_{\sigma} |\overline{E}(t_1,\sigma)|= \| \Dh(t_1,\rule{0.2em}{0.2em})- D(t_1,\rule{0.2em}{0.2em}) \|_{L^\infty}.$$ Thus, taking the max versus $s$, we get
$$
\| \Dh(t_2,\rule{0.2em}{0.2em})- D(t_2,\rule{0.2em}{0.2em}) \|_{L^\infty}\leq \| \Dh(t_1,\rule{0.2em}{0.2em})- D(t_1,\rule{0.2em}{0.2em}) \|_{L^\infty}
.
$$

Since
$$
\overline{E}(t,s) = \overline{E}(0,s) e^{-k \int_0^t \big (\overline{\Gamma}(\tau ,s) \|g_\tau (\phi(\tau,s))\|^2\big) d\tau }
$$
we have ($\phi(0,s)\equiv s$).
$$
|\overline{E}(t,s)|\leq |{E}(0,s)| e^{-k \bar\gamma \int_0^t  \|g_\tau (\phi(\tau,s))\|^2~d\tau }
.
$$
Take $p>0$. Then
$$
\int_{\SSS^2} \big|E(t,\eta)\big|^p d\sigma_\eta
=
\int_{\SSS^2} \big|\overline{E}(t,s)\big|^p \det\left(\dv{\phi}{s} (t,s)\right) d\sigma_s
.
$$
By assumption  $\dv{\phi}{s}$ is bounded.  Thus exists $C>0$ such that
$$
 \|{E}(t,\rule{.2em}{.2em})\|_{L^p}=\int_{\SSS^2} \big|E(t,\eta) \big|^p d\sigma_\eta
\leq
 C \int_{\SSS^2} \big|\overline{E}(t,s)\big|^p  d\sigma_s
.
$$

When $\int_0^{+\infty } \|g_\tau(\phi(\tau,s))\|^2 d\tau =+\infty$, for each $s$ we have  $\lim_{t\mapsto +\infty} \overline{E}(t,s)=0$. Moreover $|\overline{E}(t,s)|$ is uniformly bounded by the $L^\infty $ function $E(0,s)$. By Lebesgue dominate convergence theorem  $\lim_{t\mapsto +\infty} \|\overline{E}(t,\rule{.2em}{.2em})\|_{L^p}=0$. Previous inequality leads to
$\lim_{t\mapsto +\infty}\|{E}(t,\rule{.2em}{.2em})\|_{L^p}=0$.

When, for $t>T$,  $\int_0^{T } \|g_\tau(\phi(\tau,s))\|^2 d\tau \geq \lambda t $,
we have, for all $s\in\SSS^2$, $|\overline{E}(t,s)|\leq |E(0,s)|e^{-k\bar\gamma\lambda t}$. Thus, for all $s\in\SSS^2$ we get
$|\overline{E}(t,s)| \leq \|{E}(0,\rule{.2em}{.2em})\|_{L^\infty}$. Since $\eta\mapsto \phi^{-1}(t,\eta)$ is a diffeomorphism of $\SSS^2$, we get finally, for all $\eta\in\SSS^2$,
$|{E}(t,\eta)| \leq \|{E}(0,\rule{.2em}{.2em})\|_{L^\infty}e^{-k\bar\gamma\lambda t}$. This proves
$\|{E}(t,\rule{.2em}{.2em})\|_{L^\infty} \leq \|{E}(0,\rule{.2em}{.2em})\|_{L^\infty}e^{-k\bar\gamma\lambda t}$.
\end{proof}

\subsection{Asymptotic observer based on rough depth estimation ($\GHS$)}
\label{ssec:ObsGHS}
Instead of relying the  observer on estimation $\VHS$,  we can base it on $\GHS$. Then~\eqref{eq:obserIntrinsic} becomes
($k>0$)
\begin{equation}\label{eq:obserIntrinsic_Bis}
\dv{ \Dh}{t} =- \nabla \Dh \cdot(f_t+\GHS g_t)- v\cdot\eta + k  ( 1 -\Dh \GHS)
\end{equation}
that reads in  pinhole coordinates
\begin{multline}
\dv{\Dh}{t}=-\dv{\Dh}{z_1}(f_1+\GHS g_1)-\dv{\Dh}{z_2}(f_2+\GHS g_2)\\
-(z_1v_1+z_2v_2+v_3)+ k  (1-\Dh \GHS)
.
\label{eq:observerPH_Bis}
\end{multline}
For this observer we have the following convergence result.
\begin{thm}\label{thm:error_Bis}
Take assumptions of theorem~\ref{thm:error} concerning the scene surface $\Sigma$,  $\Gamma=1/D$, $v$ and  $\omega$.
Consider the observer~\eqref{eq:obserIntrinsic_Bis} where $\GHS$ coincides with $\Gamma$ and  where  the  initial condition   is $C^1$ versus $\eta$. Then
 $\forall t\geq 0$, the solution $\Dh(t,\eta)$ of~\eqref{eq:obserIntrinsic_Bis}  exists, is unique, remains $C^1$ versus $\eta$ and
    $$
    \| \Dh(t,\rule{0.2em}{0.2em})- D(t,\rule{0.2em}{0.2em}) \|_{L^\infty} \leq
    e^{-k\bar\gamma t}
    \| \Dh(0,\rule{0.2em}{0.2em})- D(0,\rule{0.2em}{0.2em}) \|_{L^\infty}
    $$
(exponential convergence in $L^\infty$ topology)

\end{thm}

The proof, similar to the one of~theorem~\ref{thm:error}, is  left to the reader.

\section{Simulations and numerical implementations}
\label{sec:implement}
\subsection{Sequence of synthetic images and method of comparison}
\label{ssec:sequence}
\begin{figure}[thpb]
      \centering
      \includegraphics[width=0.7\textwidth,height=0.7\textwidth]{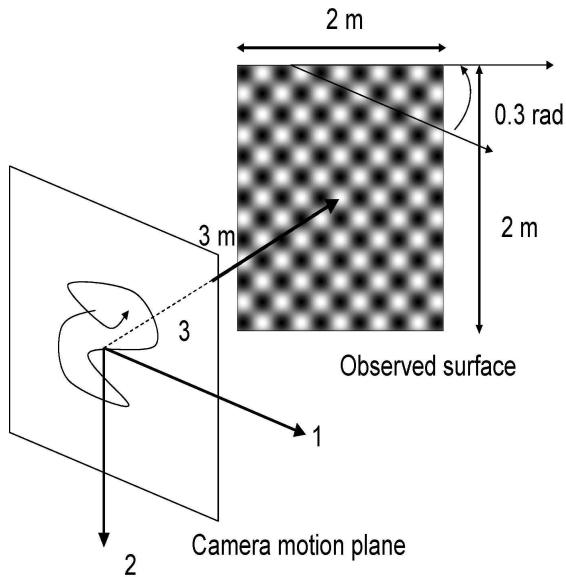}
      \caption{Virtual setup used to generate the sequence of images processed in \ref{sec:implement}}
      \label{fig:synthetic_film_setup}
\end{figure}

The non-linear asymptotic observers described in section \ref{sec:Observer} are tested on a sequence of synthetic images characterized by the following:
\begin{itemize}
\item \textit{virtual camera}: the size of each image is 640 by 480 pixels, the frame rate of the sequence is 60 Hz and the field of view is 50 deg by 40 deg;
\item \textit{motion of the virtual camera}: it consists of two combined translations in a vertical plane ($v_3=\omega_1=\omega_2=\omega_3=0$), and the velocity profiles are sinusoids with magnitude 1~$m.s^{-1}$, and different pulsations ($\pi$ for $v_1$ and $3\pi$ for $v_2$);
\item \textit{virtual scene}: it consists of a 4 $m^2$-plane placed at 3~m and tipped of an angle of 0.3 rad with respect to the plane of camera motion; the observed plane is virtually painted with a gray pattern, whose intensity varies in horizontal and vertical directions as a sinusoid function;
\item \textit{generation of the images}: each pixel of an image has an integer value varying from $1$ to $256$, directly depending on the intensity of the observed surface in the direction indexed by the pixel, to which a normally distributed noise varying with mean 0 and standard deviation $\sigma$ is added.
\end{itemize}
The virtual setup used to generate the sequence of images is represented in Fig.\ref{fig:synthetic_film_setup}.

To compare the performances of both methods, we use the global error rate in the estimation of $D$, defined as
\begin{equation}
E=\int_{ \mathcal I}\tfrac{\abs{\Dh(t,\eta)-D(t,\eta)}}{D(t,\eta)} d\sigma_\eta~/ \int_{ \mathcal I} d\sigma_\eta
\label{eq:Error}
\end{equation}
where $D$ is the true value of the depth field, $\Dh$ is the estimation computed by any of the proposed methods and $\mathcal I$ is the image frame.

\subsection{Implementation of the depth estimation based on optical flow measures ($\VHS$)}
\label{ssec:Obs1}
\begin{figure}[thpb]
      \centering
      \includegraphics[width=\textwidth]{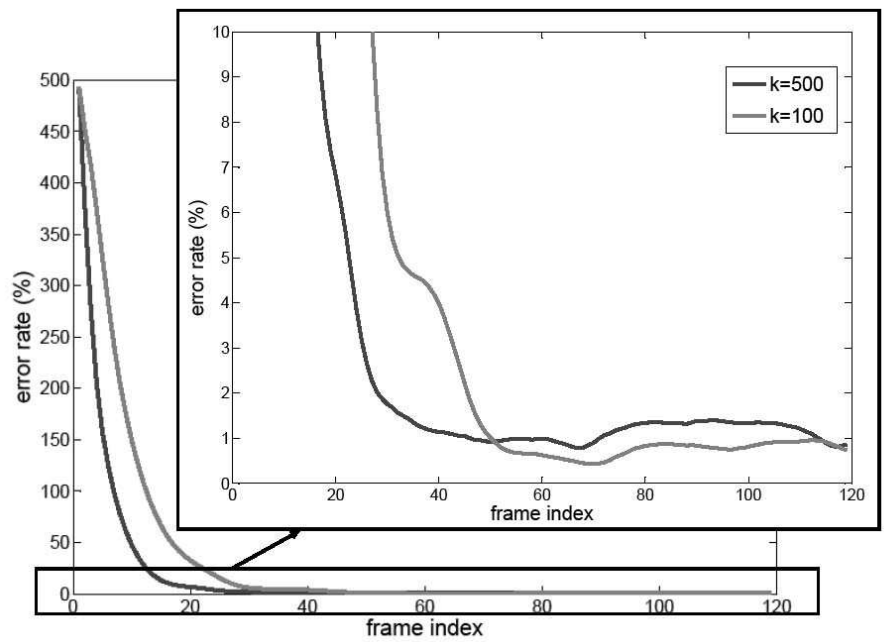}
      \caption{Relative estimation errors of the depth field $D$ estimated by the asymptotic observer \eqref{eq:observerPH} filtering the optical flow input $\VHS$ obtained by Horn-Schunck method, for different correction gains k. The noise corrupting the image data is normally distributed, with mean $\mu=0$ and standard deviation $\sigma=1$.}
      \label{fig:errObs1sigma1total}
\end{figure}

\begin{figure}[thpb]
      \centering
      \includegraphics[width=\textwidth]{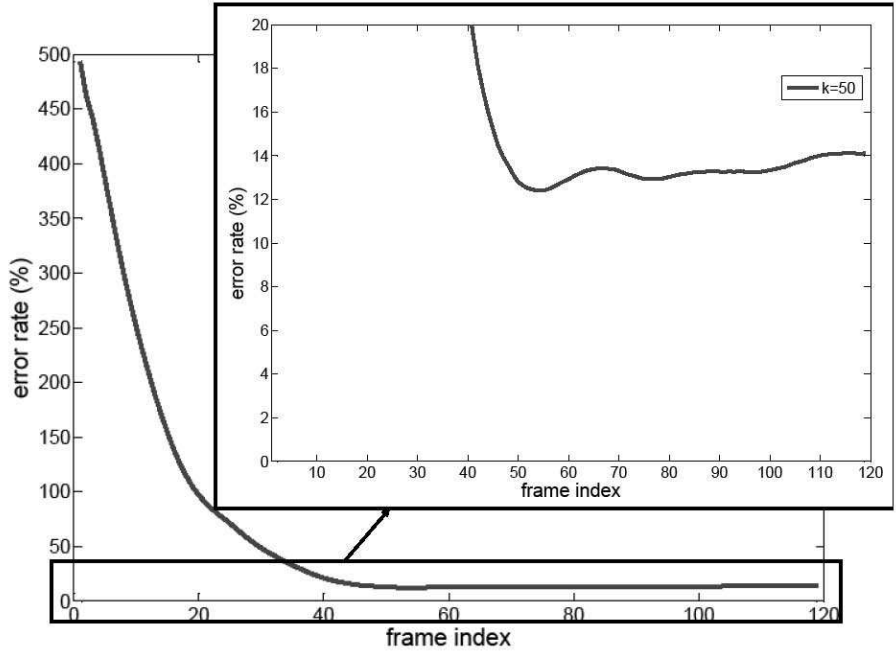}
      \caption{Relative estimation errors of the depth field $D$ estimated by the asymptotic observer \eqref{eq:observerPH} filtering the optical flow input $\VHS$ obtained by Horn-Schunck method.The noise corrupting the image data is normally distributed, with mean $\mu=0$ and standard deviation $\sigma=20$.}
      \label{fig:errObs1sigma20total}
\end{figure}

We test on the sequence described in \ref{ssec:sequence} the depth estimation characterized by the partial differential equation \eqref{eq:intrinsicEDP}.
 The optical flow input $\VHS$ ($V_1$ and $V_2$ components) is computed by a classical Horn-Schunck method. Note that  convergence theorem~\ref{thm:error} assumes that the domain of definition of the image was the entire unit sphere $\SSS^2$. Here the field of view of our virtual camera limits this domain to a portion of the sphere $\mathcal I$. However, the motion of our virtual camera ensures that most of the points of the scene appearing in the first image stay in the field of view of the camera during the whole sequence. The convergence of the method is only ensured for these points, and Neumann boundary conditions are chosen at the borders where optical flow points toward the inside of the image:
\begin{equation}
\dv{\Dh}{n}=0 \text{ if } n \cdot \VHS <0
\end{equation}
The observer gain $k>0$ is chosen in accordance with scaling considerations. Setting $k=500$ $s.m^{-2}$ provides a rapid convergence rate: we see on Fig.~\ref{fig:errObs1sigma1total} that after a few  frames, the initial relative error (blue curve) is reduced by $1/3$. Setting $k=50$  $s.m^{-2}$ is more reasonable when dealing with noisy data: on Fig.~\ref{fig:errObs1sigma20total} initial relative error is reduced by $1/3$ after around $20$ frames.

More precisely, the standard deviation $\sigma$ of the noise added to the synthetic sequence of images is 1. The gains $k=500$ and $k=100$ are successively tested, and the associated error rates for $\Dh$ are plotted in Fig. \ref{fig:errObs1sigma1total}. As expected, the convergence is more rapid for a larger gain but at convergence, i.e., after 40 frames,  the  relative  errors  are similar and below 1.5\%.

To test robustness  when dealing with noisy data, the standard deviation $\sigma$ is magnified by 20. The correction gain is tuned to $k=50$ $s.m^{-2}$. The converged errors after $40$ frames  significantly increases and yet stays  between 12 and 14 \%. Note that such permanent errors can not decrease since such noise  level first affects the optical flow estimation  $\VHS$ that feeds the asymptotic  observer. Compared to its true value $V$, the error level  for  $\VHS$ is about 15 \%. These results underline the fact that this approach is sensitive to input optical flow measures, but not directly to noise corrupting the image data.

\subsection{Implementation of the asymptotic observer based on rough depth estimation ($\GHS$)}
\label{ssec:Obs2}

Subsequently, the observer described by \eqref{eq:observerPH_Bis} was applied to the same sequence. The input depth field $\GHS$ is obtained as the output of  \eqref{eq:intrinsicEDP}.
To adapt the numerical method to this model, we make the small-angle approximation by neglecting the second order terms $z$ (we neglect the curvature of $\SSS^2$ and consider that the camera image corresponds to a small part of $\SSS^2$ that can be approximated by a small  Euclidean rectangle; the error  of this approximation is smaller than 3\% for such rectangular image):
 \eqref{eq:intrinsicEDP} becomes
\begin{equation}
G^2\Gamma+FG=\alpha^{2}\left( \dvv{\Gamma}{z_1} + \dvv{\Gamma}{z_2}  \right)
\label{eq:CVb}
\end{equation}
$F$ and $G$ are computed using angular and linear velocities, $\omega$ and $v$, and differentiation (Sobel)  filters directly applied on the image data $y(t,z_1,z_2)$ . $\Delta \Gamma$ is approximated by  the difference between the weighted mean $\bar \Gamma$ of $\Gamma$  on the neighboring pixels and its value at the current pixel. The resulting linear system in $\Gamma$ is solved by  the Jacobi iterative scheme, with an initialization provided by the previous estimation. The regularization parameter $\alpha$ is chosen accordingly to scaling considerations and taking into account the magnitude of noise: $\alpha=300$ $m.s^{-1}$ provides a convergence in about 5 or 6 frames for relatively clean image data.
As for the observer \eqref{eq:observerPH_Bis}, the correction gain $k=50$ $s.m^{-1}$ enables a convergence in around 20 frames.

\begin{figure}[thpb]
      \centering
      \includegraphics[width=\textwidth]{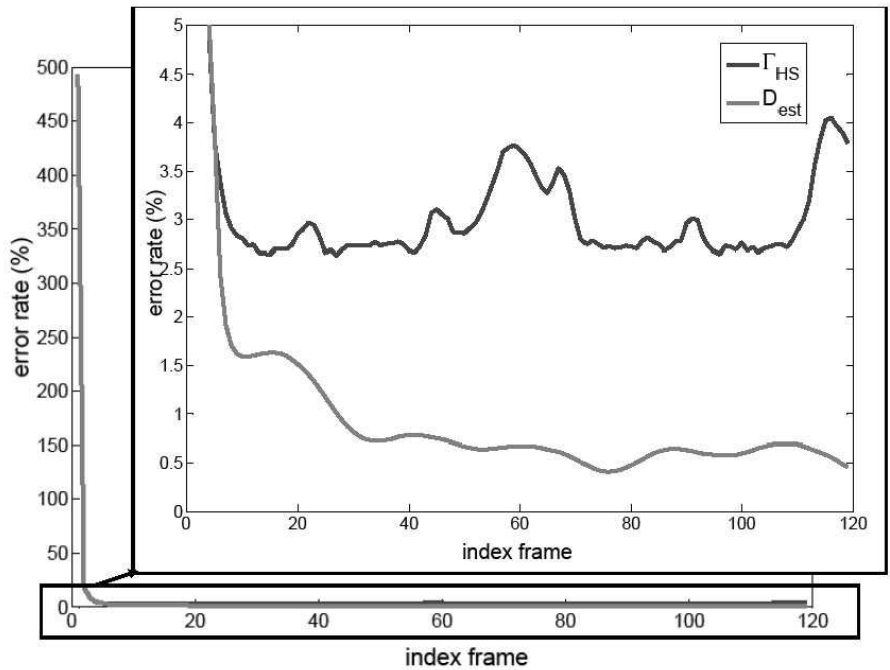}
      \caption{Relative estimation error of $\GHS$ (blue), using the depth estimation inspired by Horn-Schunck, described in \ref{ssec:HSDepth} and of $D$ (red) estimated by the asymptotic observer \eqref{eq:observerPH_Bis} filtering $\GHS$.The noise corrupting the image data is normally distributed, with mean $\mu=0$ and standard deviation $\sigma=1$.}
      \label{fig:errObs2sigma1total}
\end{figure}

As in section \ref{ssec:Obs1}, we test the observer \eqref{eq:observerPH_Bis} for different levels of noise corrupting the image data. For $\sigma=1$, the error rates associated to the input depth $\GHS$ and to the estimated depth $\Dh$ are plotted in Fig. \ref{fig:errObs2sigma1total}.
After only 6 images of the sequence, the error rate for $\GHS$ is smaller than 4\% and stays below this upper bound for the rest of the sequence. On the downside, the error rate stays larger than 2.5\%. On the contrary, the error rate associated to $\Dh$ keeps decreasing, and reaches the minimal value of 0.5 \%.

\begin{figure}[thpb]
      \centering
      \includegraphics[width=\textwidth]{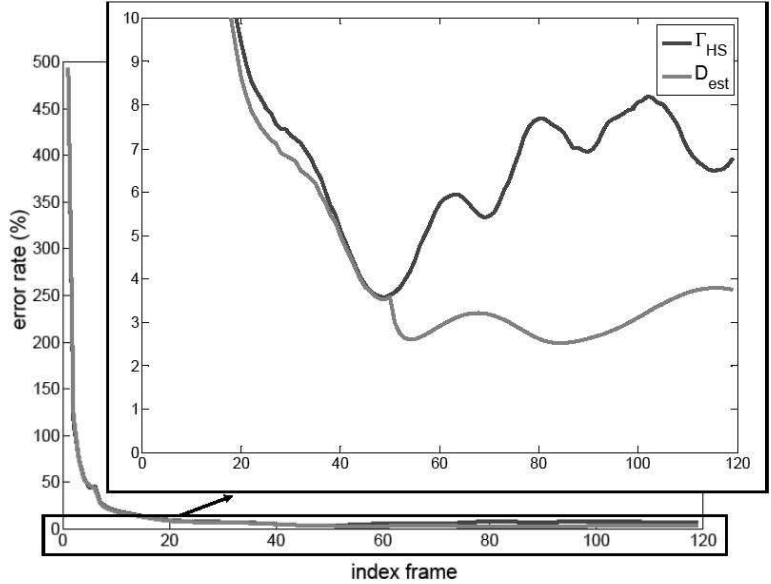}
      \caption{Relative estimation error of $\GHS$ (blue), using the depth estimation inspired by Horn-Schunck, described in \ref{ssec:HSDepth} and of $D$ (red) estimated by the asymptotic observer \eqref{eq:observerPH_Bis} filtering $\GHS$.The noise corrupting the image data is normally distributed, with mean $\mu=0$ and standard deviation $\sigma=20$.}
      \label{fig:errObs2sigma20total}
\end{figure}

For $\sigma=20$, the error rates associated to the input depth $\GHS$ and to the estimated depth $\Dh$ are plotted in Fig. \ref{fig:errObs2sigma20total}. For the computation of $\GHS$, the diffusion parameter $\alpha$ is increased to $1000$ $m.s^{-1}$ to take into account such  stronger noise.
The observer filters the error associated to $\GHS$ (between 4 and 8 \%) to provide a 3 \%  accuracy. The results show a good robustness to noise for this observer.

\subsection{Experiment on real data}
To realize the experiments, a camera was fixed on a motorized trolley traveling back and forth on a 2 meter-linear track in about 6 seconds. The resolution of the encoder of the motor enables to know the position and the speed of the trolley with a micrometric precision.  The camera is a Flea2 - Point Grey Research VGA video cameras (640 by 480 pixels) acquiring data at 20.83 fps, with a Cinegon 1.8/4.8 C-mount lens, with an angular field of view of approximately 50 by 40 deg, and oriented orthogonally to the track.
The scene is a static work environment, with desks, tables, chairs, lamps, lit up by electric light plugged on the mains, with frequency 50 Hz. The acquisition frame rate of the cameras produces an aliasing phenomenon on the video data at 4.17 Hz. In other words, the light intensity in the room is variable, at a frequency that can not be easily ignored, which does not comply with the initial hypothesis \eqref{eq:flot}. However, the impact of this temporal dependence in the equations can be reduced by a normalization of the intensity of the images such as $
y(x,y)=\frac{y(x,y)}{\bar{y}}
$
where $x$ and $y$ are the horizontal and vertical indexes of the pixels in the image, $y(x,y)$ is the intensity of this pixel and $\bar{y}$ is the mean intensity on the entire image.

\begin{figure}[thpb]
      \centering
      \includegraphics[width=\textwidth]{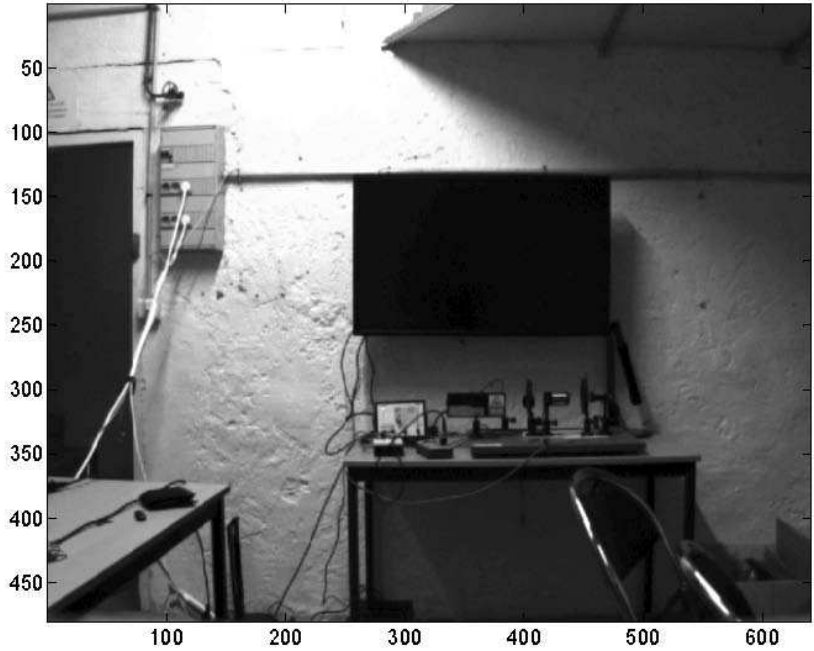}
      \caption{The static scene as seen by the camera}
      \label{fig:film}
\end{figure}

\begin{figure}[thpb]
      \centering
      \includegraphics[width=\textwidth]{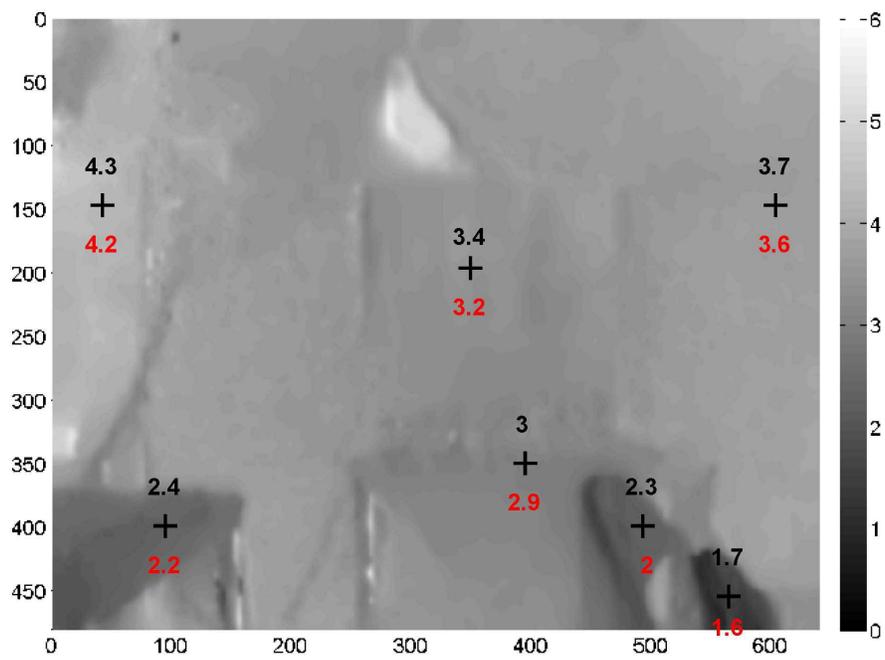}
      \caption{Estimation of the depth field associated with the image shown in Fig.\ref{fig:film}. Depth is associated to a gray level, whose scale in meters is on the right. Some estimates are extracted from the entire field (in black) and compared to real measures (in red).}
      \label{fig:depth_film}
\end{figure}

The depth field was estimated via the asymptotic observer~\eqref{eq:observerPH} based on optical flow measures. The components of the optical flow were computed using a high quality algorithm based on TV-$L^1$ method (see \cite{TVL1} for more details). The correction gain was tuned to $k=100$~$s.m^{-2}$.
An example of image data is shown in Fig.\ref{fig:film}, and the depth estimate associated to that image at the same time $t*$ is shown in Fig.\ref{fig:depth_film}. At that specific time $t^\star\approx 8$~s, the trolley already traveled once along the track and is on its way back toward its starting point. Some specific estimates are extracted from the whole depth field (two tables, two chairs, a screen, two walls) and highlighted in black; they are compared to real measures taken in the experimental room (in red):  the estimate depth profile  $\Dh(t^\star,\cdot)$ exhibits a strong correlation with these seven  punctual reference values of $D(t^\star,\cdot)$; the global appearance of the depth field looks very realistic.


\section{Conclusions and future works}

In section~\ref{sec:probstat}, we recalled  a system of partial differential equations, describing the  invariant dynamics of brightness and depth smooth fields under the assumptions of a static and Lambertian environment. We proposed in section~\ref{sec:VarMeth} an adaptation of optical flow algorithms that take the best advantage of the $SO(3)$-invariance of these equations and the knowledge of camera dynamics: it yielded an $SO(3)$-invariant  variational method to directly estimate the  depth field. In section~\ref{sec:Observer}, we proposed two asymptotic observers, respectively based on optical flow and on depth estimates. We proved their convergence under geometric  and  persistent excitation assumptions. On synthetic images, we showed in section~\ref{sec:implement} that the variational method converges rapidly, but its performance is highly dependent of the noise level whereas both asymptotic observers filter this noise. These asymptotic obersevers based on image processing of the entire field of view of the camera seems to be an interesting tool to dense range estimation and a complement to methods based on feature tracking.


\begin{thebibliography}{10}
\providecommand{\url}[1]{#1}
\csname url@samestyle\endcsname
\providecommand{\newblock}{\relax}
\providecommand{\bibinfo}[2]{#2}
\providecommand{\BIBentrySTDinterwordspacing}{\spaceskip=0pt\relax}
\providecommand{\BIBentryALTinterwordstretchfactor}{4}
\providecommand{\BIBentryALTinterwordspacing}{\spaceskip=\fontdimen2\font plus
\BIBentryALTinterwordstretchfactor\fontdimen3\font minus
  \fontdimen4\font\relax}
\providecommand{\BIBforeignlanguage}[2]{{%
\expandafter\ifx\csname l@#1\endcsname\relax
\typeout{** WARNING: IEEEtran.bst: No hyphenation pattern has been}%
\typeout{** loaded for the language `#1'. Using the pattern for}%
\typeout{** the default language instead.}%
\else
\language=\csname l@#1\endcsname
\fi
#2}}
\providecommand{\BIBdecl}{\relax}
\BIBdecl

\bibitem{Smith90}
R.~Smith, M.~Self, and P.~Cheeseman, \emph{Estimating uncertain spatial
  relationships in robotics}.\hskip 1em plus 0.5em minus 0.4em\relax New York,
  NY, USA: Springer-Verlag New York, Inc., 1990, pp. 167--193.

\bibitem{EKFMono}
J.~Civera, A.~Davison, and J.~Montiel, ``Inverse depth parametrization for
  monocular slam,'' \emph{Robotics, IEEE Transactions on}, 2008.

\bibitem{Montemerlo03a}
M.~Montemerlo, S.~Thrun, D.~Koller, and B.~Wegbreit, ``{FastSLAM} 2.0: An
  improved particle filtering algorithm for simultaneous localization and
  mapping that provably converges,'' in \emph{
  International Joint Conference on Artificial Intelligence IJCAI}, 2003.

\bibitem{StrasdatMD10}
H.~Strasdat, J.~M.~M. Montiel, and A.~J. Davison, ``Real-time monocular slam:
  Why filter?'' in \emph{ICRA}, 2010, pp. 2657--2664.

\bibitem{Strasdat-RSS-10}
H.~Strasdat, J.~M.~M. Montiel, and A.~Davison, ``Scale drift-aware large scale
  monocular slam,'' in \emph{Proceedings of Robotics: Science and Systems},
  Zaragoza, Spain, June 2010.

\bibitem{AbdursulIG04}
R.~Abdursul, H.~Inaba, and B.~K. Ghosh, ``Nonlinear observers for perspective
  time-invariant linear systems,'' \emph{Automatica}.

\bibitem{DahlWLH10}
O.~Dahl, Y.~Wang, A.~F. Lynch, and A.~Heyden, ``Observer forms for perspective
  systems,'' \emph{Automatica}.

\bibitem{Heyden09}
A.~Heyden and O.~Dahl, ``Provably convergent on-line structure and motion
  estimation for perspective systems,'' in \emph{Computer Vision Workshops,
  2009 IEEE 12th International Conference on}, 2009.

\bibitem{chen02}
X.~Chen and H.~Kano, ``A new state observer for perspective systems,''
  \emph{Automatic Control, IEEE Transactions on}, 2002.

\bibitem{Dixon03}
W.~Dixon, Y.~Fang, D.~Dawson, and T.~Flynn, ``Range identification for
  perspective vision systems,'' \emph{Automatic Control, IEEE Transactions on},
  vol.~48, no.~12, pp. 2232 -- 2238, dec. 2003.

\bibitem{Karagiannis05}
D.~Karagiannis and A.~Astolfi, ``A new solution to the problem of range
  identification in perspective vision systems,'' \emph{Automatic Control, IEEE
  Transactions on}, vol.~50, no.~12, pp. 2074 -- 2077, dec. 2005.

\bibitem{deluca08}
A.~D. Luca, G.~Oriolo, and P.~R. Giordano, ``Feature depth observation for
  image-based visual servoing: Theory and experiments.'' \emph{I. J. Robotic
  Res.}, pp. 1093--1116, 2008.

\bibitem{Astolfi10}
M.~Sassano, D.~Carnevale, and A.~Astolfi, ``Observer design for range and
  orientation identification,'' \emph{Automatica}.

\bibitem{MatthiesKS89}
L.~Matthies, T.~Kanade, and R.~Szeliski, ``Kalman filter-based algorithms for
  estimating depth from image sequences,'' \emph{International Journal of
  Computer Vision}, vol.~3, no.~3, pp. 209--238, 1989.

\bibitem{Hoilund10}
C.~Hoilund, T.~Moeslund, C.~Madsen, and M.~Trivedi, ``Improving stereo camera
  depth measurements and benefiting from intermediate results,'' in \emph{IEEE
  Intelligent Vehicles Symposium}, 2010, pp. 935 --940.

\bibitem{Herisse10}
B.~Herisse, S.~Oustrieres, T.~Hamel, R.~E. Mahony, and F.-X. Russotto, ``A
  general optical flow based terrain-following strategy for a vtol uav using
  multiple views,'' in \emph{ICRA}, 2010, pp. 3341--3348.

\bibitem{Bonnabel-Rouchon2009}
S.~Bonnabel and P.~Rouchon, ``Fusion of inertial and visual : a geometrical
  observer-based approach,'' \emph{2nd Mediterranean Conference on Intelligent
  Systems and Automation (CISA'09)}, 2009.

\bibitem{Murray09}
A.~Censi, S.~Han, S.~B. Fuller, and R.~M. Murray, ``A bio-plausible design for
  visual attitude stabilization,'' in \emph{CDC}, 2009, pp. 3513--3520.

\bibitem{HS81}
B.~Horn and B.~Schunck, ``Determining optical flow,'' \emph{Artificial
  Intelligence}, vol.~17, pp. 185--203, 1981.

\bibitem{convHS}
A.~Mitiche and A.~reza Mansouri, ``On convergence of the horn and schunck
  optical-flow estimation method,'' \emph{IEEE Transactions on Image
  Processing}, vol.~13, pp. 848--852, 2004.

\bibitem{TVL1}
A.~Chambolle and T.~Pock, ``A first-order primal-dual algorithm for convex
  problems with applications to imaging,'' \emph{Journal of Mathematical
  Imaging and Vision}, pp. 1--26, 2010.

\end{thebibliography}

\end{document}